# Qualitative visualization of distance information

Jobst Heitzig

Institut für Mathematik, Universität Hannover
Welfengarten 1, D-30167 Hannover, Germany

heitzig@math.uni-hannover.de

October 26, 2018

**Abstract**

Different types of two- and three-dimensional representations of a finite metric space are studied that focus on the accurate representation of the linear order among the distances rather than their actual values. Lower and upper bounds for representability probabilities are produced by experiments including random generation, a rubber-band algorithm for accuracy optimization, and automatic proof generation. It is proved that both *farthest neighbour representations* and *cluster tree representations* always exist in the plane. Moreover, a measure of *order accuracy* is introduced, and some lower bound on the possible accuracy is proved using some clustering method and a result on maximal cuts in graphs.

## 1 Introduction

The question of how distance information might be visualized is of importance for many sciences including physics, medicine, sociology, and others. Mathematicians have early studied the possibility of *embedding* a finite metric space $\underline{X}$ into other, in some sense better spaces like the Euclidean plane or 3-space. Beginning with Menger [Men28], who gave the precise criteria for $\underline{X}$ to be *isometrically* embeddable (that is, under exact preservation of the distances) into some Euclidean space, most of them have focused on mappings that map $\underline{X}$ into some standard space in a "quantitative" manner. The goal in this field of research, known under the name *metric scaling,* is to preserve the *values* of the distances as good as possible, that is, to minimize a certain error, known as "stress" (cf. [She62]).

The aim of this paper is to study more "qualitative" kinds of visualization of distance data. In contrast to metric scaling, we will not be interested in the actual values of distances but rather in their comparison. Considering only the linear *order* among the distances instead of their value, a measure of *order accuracy* of a representation is introduced. Unlike stress, order accuracy has an easy interpretation as a certain probability of correctness. After an experimental exploration of different types of representations, a lower bound on the possible accuracy of plane representations will be proved using some clustering method and a result on maximal cuts in graphs. The experimental methods include random generation, optimization of accuracy by a rubber-band algorithm, and automatic proof generation. All results are summarized in Table 1.

## 2   Order accuracy

Throughout this paper, $\underline{X} = (X, d)$ is a *finite metric space,* that is, $X$ is finite, and $d : X^2 \to [0, \infty]$ fulfils $d(x, y) = d(y, x), d(x, y) + d(y, z) \geqslant d(x, z)$, and $d(x, y) = 0$ if and only if $x = y$. However, one advantage of the following approach is that it also applies to any *finite, symmetric distance set* in the sense of [Hei98] and [Hei02], which is a far more general type of object than a metric space. For the sake of simplicity, we will also assume that $X$ equals the set $n = \{0, \ldots, n-1\}$ of non-negative integers, and that the pairwise distances between the points of $\underline{X}$ are all different, that is, $d(x, y) = d(x', y') > 0$ implies $\{x, y\} = \{x', y'\}$. In particular, each $x \in X$ has exactly one *nearest neighbour* $\mathrm{nn}(x) \in X$ and one *farthest neighbour* $\mathrm{fn}(x)$ which fulfil $d(x, \mathrm{nn}(x)) < d(x, y) < d(x, \mathrm{fn}(x))$ for all $y \in X \setminus \{x, \mathrm{nn}(x), \mathrm{fn}(x)\}$.

We will be mostly interested in representing the points of $X$ by points of either some Euclidean space $\mathbb{E}_m$, that is, the real vector space $\mathbb{R}^m$ with Euclidean distance, or the $L_1$-plane $\mathbb{M}_2$, that is, the set $\mathbb{R}^2$ with the "Manhattan"-distance $d(x, y) = |x_1 - y_1| + |x_2 - y_2|$.

The *order accuracy* $\alpha(f)$ of a map $f$ from $\underline{X}$ into some metric space $\underline{Y} = (Y, e)$ is defined as the probability that, of two randomly chosen pairs $\{x, y\}$ and $\{z, w\}$ of elements of $X$, the one with the larger distance in the "representation" $f$ also has the larger "original" distance. More formally,

$$\alpha(f) = \binom{\binom{n}{2}}{2}^{-1} \cdot \Big| \Big\{ \{\{x, y\}, \{z, w\}\} \subseteq \mathcal{P}(X) : \\ x \neq y,\ z \neq w,\ \{x, y\} \neq \{z, w\},\ \text{and} \\ d(x, y) < d(z, w) \iff e(fx, fy) < e(fz, fw) \Big\} \Big|.$$

Note that $2\alpha(f) - 1$ is just Kendall's rank correlation coefficient $\varrho$ between the two linear orders on the $\binom{n}{2}$ pairs $\{x, y\}$ that result when these pairs are compared with respect to either their original or their image distance. Using a variant of the merge-sort algorithm, $\varrho$ can be computed in linear-times-logarithmic time, hence $\alpha(f)$ can be computed in $O(n^2 \log n)$ time.

## 3   Order and weaker representations

An *order representation* of $\underline{X}$ in $\underline{Y}$ is some map $f : \underline{X} \to \underline{Y}$ with $\alpha(f) = 1$, that is, with $d(x, y) < d(z, w) \iff e(fx, fy) < e(fz, fw)$. Likewise, an *order representation* of a (strict) linear order $<$ on the set $\mathcal{B}(X)$ of two-element subsets of $X$ is a map $f : X \to \underline{Y}$ with $\{x, y\} < \{z, w\} \iff e(fx, fy) < e(fz, fw)$. It will be convenient to identify the metric space $\underline{X}$ with its associated linear order $<$ which is given by $\{x, y\} < \{z, w\} :\iff d(x, y) < d(z, w)$ here.

Table 1: Representable fraction of linear orders on $\mathcal{B}(n)$ for different kinds of representations and different spaces (open intervals show exact bounds, $\gtrsim$ and $\lesssim$ denote estimated bounds, the question mark denotes a conjecture)

| kind of representation | space | \multicolumn{9}{c}{$n$ (no. of points)} | | | | | | | | |
|---|---|---|---|---|---|---|---|---|---|---|
| | | 4 | 5 | 6 | 7 | 8 | 9 | 13 | 15 | any |
| order | $\mathbb{E}_2$ | **1** | $(.538, 1)$ | $\gtrsim .020$ | | | | | | |
| | $\mathbb{M}_2$ | 1 | $(.652, 1)$ | | | | | | | |
| | $\mathbb{E}_3$ | | **1** | $\gtrsim .60$ | $\gtrsim .095$ | $< 1$ | | | | |
| local order | $\mathbb{E}_2$ | | $(\mathbf{.667}, 1), \lesssim .928$ | $\lesssim .60$ | $\lesssim .21$ | $\lesssim .030$ | $\lesssim .0022$ | | | |
| | $\mathbb{M}_2$ | | $(.677, 1)$ | | | | | | | |
| extremal neighbours | $\mathbb{E}_2$ | | $(\mathbf{.883}, 1)$ | | | | | | | |
| $1^{\text{st}}$ & $2^{\text{nd}}$ nearest nbs. | $\mathbb{E}_2$ | | $(\mathbf{.933}, 1)$ | | | | | | | |
| two nearest neighbours | $\mathbb{E}_2$ | | $(\mathbf{.963}, 1)$ | | | | | | | |
| nearest neighbour | $\mathbb{E}_2$ | | | **1** | **.999**... | **.998**... | **.997**... | | | |
| nearest neighbour | $\mathbb{E}_3$ | | | | | | | **1?** | $< 1$ | |
| farthest neighbour | $\mathbb{E}_2$ | | | | | | | | | **1** |
| cluster | $\mathbb{E}_2$ | | | | | | | | | **1** |

For $\underline{Y} = \mathbb{E}_{n-2}$, there is always an order representation—there is even a map $f$ for which $e(fx, fy) = d(x, y) + C$ for some constant $C \geqslant 0$. This was proved by Cailliez [Cai83]. A random generation of five-element subsets of $\mathbb{E}_3$ confirmed this result for $n = 5$, and a similar experiment showed that all four-element metric spaces not only have an order representation in $\mathbb{E}_2$ but also in $\mathbb{M}_2$.

To get a feeling how probable a *plane* order representation is for a five-element metric space, I also repeatedly drew five-element samples from the uniform distribution on the unit square and determined the resulting order among the ten pairwise distances. In this way, of the $10! = 3\,628\,800$ linear orders on $\mathcal{B}(5)$, at least $53.8\%$ [resp. $65.2\%$] were found to have an order representation in $\mathbb{R}^2$ with the Euclidean [resp. "Manhattan"] metric. Moreover, at least $66.7\%$ [resp. $67.7\%$] had a *local order representation*, that is, a map $f : \underline{X} \to \mathbb{R}^2$ such that $\{x, y\} < \{x, z\} \iff e(fx, fy) < e(fx, fz)$ for all $x, y, z$, where again $e$ was the Euclidean [resp. "Manhattan"] metric. Judging from these empirical numbers, order representability seems to be considerably stronger than local order representability in the Euclidean case, but not in the "Manhattan" case.

Considering only the information coded in the functions nn and fn, it was also found that at least $88.3\%$ of the 10! orders had a plane *extremal neighbours representation*, that is, a map $f : \underline{X} \to \mathbb{E}_2$ such that $\text{nn}(fx) = f(\text{nn}(x))$ and $\text{fn}(fx) = f(\text{fn}(x))$ for all $x \in X$. Likewise, at least $93.3\%$ allowed for a map under which both the nearest and second-nearest neighbours were represented accurately, and another $3\%$ allowed for a map under which at least the information about which points were the two nearest to $x$ was represented accurately for all $x$ (see Table 1).

In view of the quickly growing number $\binom{n}{2}!$ of orders on $\mathcal{B}(n)$ and the limited space for storing the list of orders already found, such a random generation did not make much sense for $n > 5$. It is, however, possible to estimate some similar lower bounds at least for $n \in \{6, 7\}$ from the following experiment.

## 4 Representation by accuracy optimization

Starting with a randomly generated $f : X \to \mathbb{E}_m$, an order representation of a linear order $<$ on $\mathcal{B}(X)$ can often be produced by a stepwise maximization of order accuracy. The following optimization step proved useful: for each pair $\{x, y\}$, $\{z, w\}$ with $\{x, y\} < \{z, w\}$ and $e(fx, fy) \geqslant e(fz, fw)$, move $x, y$ towards each other by some fixed fraction of $e(fx, fy)$, and move $z, w$ away from each other by the same fixed fraction of $e(fz, fw)$. I have tested this kind of rubber-band algorithm in several ways:

(i) When $<$ was taken to be the order that corresponded to 8 or 25 independently uniformly distributed random points in the unit square, the algorithm found an order representation of $<$ in $\mathbb{E}_2$ in about $96\%$ of all cases, no matter if 8 or 25 points were taken. For 25 points, the resulting representations were almost similar to the original sets. More precisely, for each edge the quotient between its original length and its length in the representation was determined, and on average the relative difference between maximal and minimal quotient was less than $5\%$ (compared to $12\%$ for 15 points and over $60\%$ for 8 points).

(ii) When $<$ was taken from a uniform distribution of all linear orders on $\mathcal{B}(5)$, the algorithm succeeded in only $45\%$ of the cases. Since, as mentioned before, more than 53% of the orders actually have an order representation, this indicates that the algorithm is susceptible to being caught in a local optimum.

However, in both (i) and (ii), the success of the algorithm did not seem to depend on the initial state: when a cluster representation (see below) instead of a random initial state was used, only the average number of iterations needed shrinked slightly.

(iii) As in (i), but for five points in a 100-dimensional cube. Here the success rate was about 79%. Such finite subspaces of high-dimensional spaces frequently occur in multivariate statistics, for example.

(iv) Generating the orders as in (ii), an order representation in $\underline{\mathbb{E}}_3$ of six-point metric spaces was found in about 65% of 1000 cases, but of seven-point spaces in only 10.5% of 7000 cases.

The rubber-band algorithm has also been implemented as a Java applet which can be tested at

http://www-ifm.math.uni-hannover.de/~heitzig/distance.

Despite the algorithm's lack of optimality, we can use these results to estimate lower bounds for the fraction of representable orders. As the samples were large enough, one can use the approximate confidence bound that arises from the approximation of the actual binomial distribution by a normal distribution (see [Kre91]). For a sample of size $N$, $s + 1/2$ successes, and confidence niveau $\beta$, it has the form

$$\frac{s + \frac{c^2}{2} - c\sqrt{s - \frac{s^2}{N} + \frac{c^2}{4}}}{N + c^2} \quad \text{with} \quad c = \Phi^{-1}(\beta).$$

Taking $\beta = 0.995$, this leads to the following conjectured bounds:

**Conjecture 1** *In $\underline{\mathbb{E}}_3$, a six- [seven-] element metric space has an order representation with probability at least* 60% *[*9.5%*].*

For six points in $\underline{\mathbb{E}}_2$, the same method gives a conjectured lower bound of only 2% (see Table 1).

## 5  Disproving local order representability

A local order representation can also be characterized as a map that preserves the order among the three sides of any triangle. More precisely, $f : \underline{X} \to \underline{Y}$ is a local order representation if and only if for each three distinct points $x, y, z \in X$ with $d(x, y) < d(y, z) < d(z, x)$, also $e(fx, fy) < e(fy, fz) < e(fz, fx)$. Using elementary geometry, one sees that, in the Euclidean plane, the latter is equivalent to $\angle fx\,fz\,fy < \angle fy\,fx\,fz < \angle fz\,fy\,fx$ ($\star$).

Therefore, the existence of a plane local order representation for some order $<$ can be disproved by showing that a certain set of inequalities between angles in the plane has no solution. The advantage of using angles instead of distances is that the additional equations and inequalities which every $n$-point subset of the plane must fulfil are all

linear in the angles:

(i)    $\angle abc \in [0, \pi]$

(ii)    $\angle abc + \angle bca + \angle cab = \pi$

(iii)    $\angle azc \leqslant \angle azb + \angle bzc$

(iv)    $\angle azb + \angle bzc + \angle cza = 2\pi$ if $z$ is in the convex hull of $a, b, c$,

(v)    $\angle azc = \angle azb + \angle bzc$ if $b$ is "between" $a$ and $c$ as seen from $z$.

In search of a local order representation for $\underline{X}$, these linear relations together with those of type $(\star)$ enable us, starting with the largest interval $[0, \pi]$, to successively narrow down the interval of possible values of each angle. If some angle's interval becomes empty, there can be no local order representation of this order $<$. This method can also be used to disprove the existence of even weaker kinds of representations such as extremal neighbours representations.

**Example 2** Figure 1 shows a computer generated proof that the order $\{d, e\} < \{a, d\} < \cdots < \{b, d\}$ (listed on top) cannot occur among the distances between five points in the plane. Lines 1, 2, and 3 state that certain angles are smaller than $60°$, smaller than $90°$, or larger than $60°$ because they are the smallest, second smallest, or largest in their corresponding triangle, respectively. Line 4 states that only $c$ can be in the convex interior of the five points, since each of the remaining four is the farthest neighbour of some other. Lines 5–7 apply the "tripod" inequality (iii), using bounds already known from lines 1 and 2, this dependence being logged at the end of the lines. Line 8 notices a violation of (iv) so that $c$ cannot be in the convex hull of $a, b, d$. Similarly, line 9 states that also $b$ cannot be between $a$ and $d$ as seen from $c$. In line 11, (ii) is used to derive a lower bound for a second smallest angle from an upper bound for a largest angle. This is the only kind of argument the algorithm can use to derive bounds that are not just multiples of $30°$. The rest of the proof shall be clear now.

Note that the premises in lines 1–4 already follow from the information coded in the maps nn and fn alone, hence the order under consideration does not even have an extremal neighbours representation.

There is a similar example which shows that it may also be impossible in the plane to accurately represent the set of two nearest neighbours of five points. Since for disjoint five-element subsets of some metric space $\underline{X}$, the distribution of the orders that correspond to these subsets are independent, we have:

**Corollary 3** *For an $n$-element metric space, the probability of a plane extremal neighbours representation shrinks exponentially for $n \to \infty$.*

To get explicit upper bounds for local representability, I tested several thousand randomly generated orders with this algorithm. For five points, 795 out of 10 000 orders could be shown to have no plane local order representation in this way. Using again estimated confidence bounds with $\beta = .995$, this results in an estimated upper bound of .928 for the fraction of plane locally order representable orders on $\mathcal{B}(5)$. For $n = 6, 7, 8$, and 9, the corresponding numbers were 4156 out of 10 000, 3627 out of 4500, 11 690 out of 12 000, and 9990 out of 10 000, respectively, resulting in the upper bounds shown in Table 1.

Figure 1: A computer generated non-representability proof.

```
TEST OF EDGE ORDER de < ad < ac < ab < ce < be < bc < cd < ae < bd
USING ONLY EXTREMAL NEIGHBOURS INFORMATION

legend: points are labeled a,b,c,d,e
        xy is a segment, xyz is a triangle, x:yz is the angle in xyz at vertex x
        x:ywz means that x:yz=x:yw+x:wz
                                                                        follows
 line type       proposition                                            from
--------------------------------------------------------------------------------
   1. smallest   a:de,b:ad,b:de,c:ad,c:de,d:bc,e:ab,e:ac < 60
   2. dominated  a:be,a:ce,b:ac,b:cd,c:ab,d:ab,d:ac,d:be,d:ce,e:ad < 90
   3. largest    a:bc,a:bd,a:cd,b:ae,c:ae,c:bd,d:ae,e:bd,e:cd > 60
   4. on bndry   a,b,d,e since in fn[X]
   5. tripod     a:bd <=a:be+a:de < 90+60= 150                          2.1.
   6. tripod     a:cd <=a:ce+a:de < 90+60= 150                          2.1.
   7. tripod     b:ae <=b:ad+b:de < 60+60= 120                          1.1.
   8.        not c in abd since c:ad+c:bd+c:ab<360                      1.0.2.
   9.        not c:abd since c:ad<c:ab+c:bd                             1.0.3.
  10. tripod     c:ae <=c:ad+c:de < 60+60= 120                          1.1.
  11. larger     a:be > (180-b:ae)/2>(180-120)/2= 30                    7.
  12. larger     a:ce > (180-c:ae)/2>(180-120)/2= 30                    10.
  13.        not a:cbe since a:ce<a:bc+a:be                             2.3.11.
  14.        not a:bce since a:be<a:ce+a:bc                             2.12.3.
  15.      hence a:bec                                                  4.13.14.

CASE ANALYSIS using points a,bcd:

  16. (i) ASSUMING a:bcd...
  17. sum        a:bd =a:bc+a:cd > 60+60= 120                           16.3.3.
  18. sum        a:bc =a:bd-a:cd < 150-60= 90                           16.5.3.
  19. tripod     a:be >=a:bd-a:de > 120-60= 60                          17.1.
  20.        not a:bec since a:bc<a:ce+a:be                             18.12.19.
  21.      hence a in bce                                               14.13.20.
  22. contradiction!                                                    21.4.
  23. (ii) ASSUMING a:cbd...
  24. sum        a:cd =a:bc+a:bd > 60+60= 120                           23.3.3.
  25. sum        a:bc =a:cd-a:bd < 150-60= 90                           23.6.3.
  26. tripod     a:ce >=a:cd-a:de > 120-60= 60                          24.1.
  27.        not a:bec since a:bc<a:be+a:ce                             25.11.26.
  28.      hence a in bce                                               13.14.27.
  29. contradiction!                                                    28.4.
  30. (iii) ASSUMING a:bdc...
  31.        not d:acb since a:bdc                                      30.
  32.        not d:abc since a:bdc                                      30.
  33.      hence d:bac                                                  31.4.32.
  34.        not c:bad since a:bdc                                      30.
  35.      hence c:adb                                                  8.9.34.
  36. new sum    c:abd since ad diag in cabd                            30.33.
  37. new circ   d in abc since a:bdc and c:adb                         30.35.
  38. contradiction!                                                    37.4.
  39. (iv) ASSUMING a in bcd...
  40. contradiction!                                                    39.4.

CONTRADICTION in all four cases!
```

Figure 2: A "universal" nearest neighbour graph of nine points in the plane

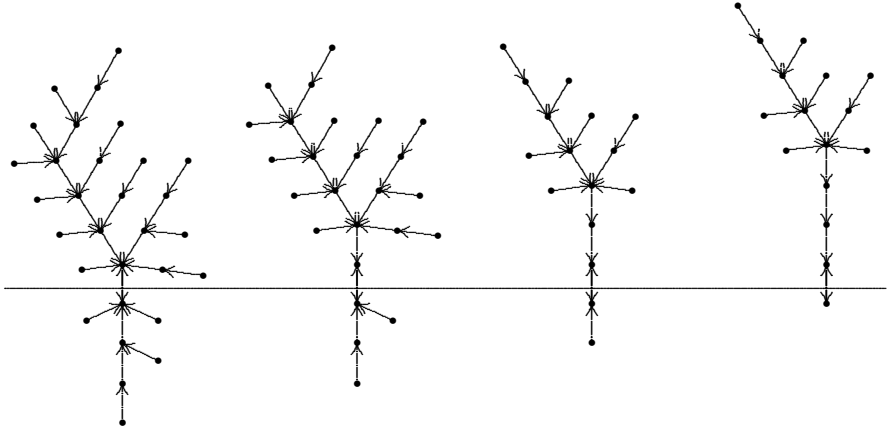

**Conjecture 4** *In $\underline{\mathbb{E}}_2$, a six-element metric space has a local order representation with probability at most* $60\%$.

This fast vanishing of the probability of plane local order representability on the one hand shows that the above algorithm is quite successful, and on the other hand motivates the study of even weaker kinds of plane representation.

# 6 Nearest and farthest neighbour representations

The directed graph $G_{\mathrm{nn}}(X)$ with vertex set $V(G) = X$ and edge set $E(G) = \{(x, \mathrm{nn}(x)) : x \in X\}$ is known as the *nearest neighbour graph* of $\underline{X}$. Asymptotic properties of nearest neighbour graphs of subsets of the plane have been studied in [EPY97]. The *farthest neighbour graph* of $\underline{X}$ is defined similarly. By a *down-tree* I mean a finite connected digraph all of whose vertices have out-degree one, except for a root vertex with out-degree zero.

**Proposition 5** *A finite digraph $G$ is a nearest [farthest] neighbour graph of a metric space if and only if each of its components is a disjoint union of two down-trees whose roots are joined by a double edge.*

Since the *proof* is easy but quite technical, it is omitted here.

The digraphs characterized by this result will be called *bi-rooted forests* in the sequel, and a pair of roots will be called a *bi-root* for short. A *proper child* of a vertex $x$ in a digraph is a vertex $y$ for which there is an edge $(y, x)$ but no edge $(x, y)$.

**Proposition 6** *A bi-rooted forest of size at most nine occurs as a nearest neighbour graph in the plane if and only if no vertex has more than four proper children.*

*Proof.* Let $G$ be a bi-rooted forest with $|V(G)| \leqslant 9$. If some vertex $x$ has five proper children $x_1, \ldots, x_5$, there is no nearest neighbour representation in $\mathbb{E}_2$. Otherwise, for $i \neq j$, the longest side of the triangle $x_i x_j x$ would be $x_i x_j$, hence the angle between the segments $x_i x$ and $x_j x$ would be larger than $\pi/3$. Likewise, the longest side of the triangle $x_i x \, \mathrm{nn}(x)$ is $x_i \, \mathrm{nn}(x)$, hence the angle $\angle x_i x \, \mathrm{nn}(x)$ would also be larger than $\pi/3$ which is impossible in the plane.

On the other hand, one can verify that all bi-rooted forests with at most nine vertices and without vertices that have more than four proper children fit into the "universal" forest sketched in Figure 2. Each of its four components is constructed from its two roots (joined by a double edge of length 100) by successively adding children, where the edges originating from children of order $n$ have length $100 + n$ and share a mutual angle of $(65 + i - n)°$ if they are neighboured. Since in that figure, each edge points towards the nearest neighbour, the proposition is proved. □

Using this result, it was possible to calculate the fractions of linear orders on $\mathcal{B}(n)$ with a plane nearest neighbour representation shown in Table 1. Note that for $n = 10$, the analogue of the above proposition is false, a counter-example being the bi-rooted forest consisting of two connected roots with four children each.

As for nearest neighbour representations in $\mathbb{E}_3$, it was proved by Fejes Tóth [FT43] that of $n$ points on a unit sphere in $\mathbb{E}_3$, at least two must have a distance of at most

$$\delta_n := \sqrt{4 - \operatorname{cosec}^2 \frac{n}{n-2} \frac{\pi}{6}}.$$

In particular, $\delta_{14} \approx 0.98$, hence there exist no fourteen points on the unit sphere with pairwise distance larger than one. In other words, of fourteen rays in $\mathbb{E}_3$ with a common source, at least two have an angle of at most $60°$. Therefore, a bi-rooted forest with a root that has thirteen children cannot have a representation in $\mathbb{E}_3$. In particular, not all linear orders on $\mathcal{B}(15)$ have a nearest neighbour representation in $\mathbb{E}_3$. However, one may hope that at least all linear orders on $\mathcal{B}(13)$ have a representation since there exist twelve such points on the sphere.

**Conjecture 7** *Every metric space of up to thirteen elements has a nearest neighbour representation in $\mathbb{E}_3$.*

Note that $\delta_{13} \approx 1.014 > 1$, and the empirically supported conjecture that there are no thirteen such points is still unproved—this might show that questions of representability of larger sets might also be quite difficult.

Surprisingly, a small degree at all vertices of the nearest neighbour graph does not assure plane nearest neighbour representability: Eppstein, Paterson, and Yao [EPY97] could show that for a subset $X$ of $\mathbb{E}_2$, $|X| = O(D(G_{\mathrm{nn}}(X))^5)$, where $D(G)$ is the *depth* of $G$, that is, the maximal length of a path from a vertex to the corresponding root. Using their exact bounds, one can show that for example the complete binary bi-rooted tree with $2^{66} - 2 \approx 10^{20}$ vertices does not have a nearest neighbour representation in $\mathbb{E}_2$. However, it seems likely that already far smaller binary trees fail to have one.

Eppstein et al. also showed that the expected number of components of $G_{\mathrm{nn}}(X)$ is asymptotic to approximately $0.31|X|$ if the points of $X$ are independently uniformly distributed in the unit square. More precisely, the probability for a vertex to belong to a bi-root is $6\pi/(8\pi + 3\sqrt{3}) \approx 0.6215$ in that case. From this it is also clear that the expected fraction of elements of $X$ that are not the nearest neighbour of some other element is at most $0.2785$. However, the smallest exact upper bound to this fraction is far larger:

**Proposition 8** *In any finite subset of $\mathbb{E}_2$, at most $7/9$ of its elements are not a nearest neighbour of some other element, and this bound is sharp.*

*Proof.* It is quite easy to see that the bi-rooted forest consisting of a root with four and another with three children has a nearest neighbour representation in $\mathbb{E}_2$, hence $7/9$ is possible.

On the other hand, let $C$ be a component of the nearest neighbour graph of a finite subset of the plane. Then its roots $r$ and $q$ together have $k \leqslant 7$ children, and $C$ can be constructed from these $k + 2$ vertices by subsequently adding $k_i \leqslant 4$ children to some end vertex, thereby increasing the number of end vertices by $k_i - 1$ in step $i$. Thus, the final fraction of end vertices in $C$ is

$$\frac{k + \sum_i (k_i - 1)}{(k+2) + \sum_i k_i} \leqslant \frac{7}{9}$$

since $7(k+2+\sum_i k_i) - 9(k+\sum_i(k_i-1)) = 14 - 2k + 9s - 2\sum_i k_i \geqslant 9s - 2\cdot 4s \geqslant 0$, where $s$ is the number of steps needed. □

In view of these facts about nearest neighbour graphs, the following might be a bit surprising:

**Theorem 9** *Every finite metric space has a farthest neighbour representation in $\mathbb{E}_2$.*

*Proof.* Let $G := G_{\text{fn}}(X)$ be the corresponding farthest neighbour graph, $D$ its depth, and define an infinite bi-rooted forest $H$ as follows. The vertices of $H$ are labelled $a_{jt}$ and $b_{jt}$, where $j$ is a non-negative integer and $t$ runs over all tuples of at most $D$ non-negative integers, including the empty tuple $\emptyset$. The bi-roots are the pairs $\{a_{j\emptyset}, b_{j\emptyset}\}$ with non-negative integer $j$, each vertex $a_{j(\ldots,k,m)}$ is a child of $a_{j(\ldots,k)}$, and each vertex $b_{j(\ldots,k,m)}$ is a child of $b_{j(\ldots,k)}$. In other words, $H$ has countably many isomorphic components (numbered by $j$), and each vertex has countably many children, up to depth $D$. This digraph $H$ contains an isomorphic copy of $G$, hence it suffices to give a representation of $H$. To address points of the plane, it will be convenient to identify $\mathbb{R}^2$ with the set $\mathbb{C}$ of complex numbers in the usual way.

For each non-negative integer $j$, let $C_{j0}$ and $C_{j1}$ be the circles of radius 2 with centres $c_{j0} := e^{2^{-j-1}\pi i}$ and $c_{j1} := e^{(1+2^{-j-1})\pi i}$, respectively. These curves can be parametrized using the following functions, where the coefficients $\lambda_j > 0$ will be determined later:

$$f_{j0}(\xi) := c_{j0} + 2e^{(2^{-j-1}+\lambda_j \xi)\pi} \quad \text{and} \quad f_{j1}(\xi) := c_{j1} + 2e^{(1+2^{-j-1}+\lambda_j \xi)\pi}.$$

In particular, $f_{j0}(0) = 3c_{j0}$, $f_{j1}(0) = 3c_{j1}$, $F_{j0} := f_{j0}[I] \subseteq C_{j0}$, and $F_{j1} := f_{j1}[I] \subseteq C_{j1}$, where $I = [-2^D, 2^D] \subseteq \mathbb{R}$. Now the coefficients $\lambda_j$ are chosen small enough so that $2^D \lambda_j < \pi/2$ and so that the smallest distance between the sets $F_{j0}$ and $F_{j1}$ is still larger than the largest distance between a point in $F_{j0} \cup F_{j1}$ and a point in $F_{k0} \cup F_{k1}$ for any $k \neq j$. This ensures that, for $q \in \{0, 1\}$ and all $\xi \in I$, the unique point in $\bigcup_k F_{k0} \cup F_{k1}$ that is farthest away from the point $f_{jq}(\xi)$ is the point $f_{j,1-q}(\xi/2)$. More generally, given $q \in \{0, 1\}$ and $\xi, \beta, \gamma \in I$, we have

$$|f_{jq}(\xi) - f_{j,1-q}(\beta)| > |f_{jq}(\xi) - f_{j,1-q}(\gamma)| \iff |\beta - \xi/2| < |\gamma - \xi/2| \quad (\star).$$

Using this equivalence, one sees that the following recursive definition results in a farthest neighbour representation $f$ of $H$:

$$f(a_{jt}) := f_{j,q(t)}(\xi(t)) \quad \text{and} \quad f(b_{jt}) := f_{j,1-q(t)}(-\xi(t)),$$

where the bi-roots have $q(\emptyset) := 0$ and $\xi(\emptyset) := 0$, their children have $q((m)) := 1$ and $\xi((m)) := 1 + 2^{-m}$, and all others have $q((\ldots, k, m)) := 1 - q((\ldots, k))$ and

$$\begin{aligned}\xi((\ldots, k, m)) &:= 2\xi((\ldots, k)) - (1 - 2^{-m})\big(\xi((\ldots, k)) - \xi((\ldots, k+1))\big) \\ &= (1 + 2^{-m})\xi((\ldots, k)) + (1 - 2^{-m})\xi((\ldots, k+1)).\end{aligned}$$

Because of $(\star)$, we need only verify that (i) $|0 - \xi((m))/2| < |\xi((k, \ell)) - \xi((m))/2|$, which is true because of $\xi((m)) < 2 < \xi((k, \ell))$, and that (ii)

$$|2\xi((\ldots, k)) - \xi((\ldots, k, m))| < |2\xi((\ldots, k \pm 1)) - \xi((\ldots, k, m))|,$$

where the left hand side equals $(1 - 2^{-m})c$ with $c = \big(\xi((\ldots, k)) - \xi((\ldots, k+1))\big)$, and the right hand side is the absolute value of $c + 2\big(\xi((\ldots, k \pm 1)) - \xi((\ldots, k, m))\big)$ which is larger than $c$ in the "$-$" case and smaller than $-c$ in the "$+$" case. $\square$

# 7 Cluster representations, and lower bounds for accuracy

A important question in applications of finite metric spaces is that of clustering the elements into homogeneous, mutually heterogeneous groups. Formally, a hierarchical clustering of $\underline{X}$ produces what I will call a *cluster tree* here, which can be formalized as a chain of partitions $\mathcal{P}_1, \ldots, \mathcal{P}_n$ on $X$, where $\mathcal{P}_1 = \{\{x\} : x \in X\}$ is the discrete and $\mathcal{P}_n = \{X\}$ the indiscrete partition, and each $\mathcal{P}_{k+1}$ with $k < n$ arises from $\mathcal{P}_k$ by joining two clusters, that is, replacing some $A, B \in \mathcal{P}_k$ by their union $A \cup B$. Most common clustering methods fulfil the following property $(\star)$: if $k < n$, $A, B \in \mathcal{P}_k$, $A \neq B$, and for all $a \in A$, $b \in B$, and $x, y \in X$, either $x, y \in A \cup B$, or $x, y \in C$ for some $C \in \mathcal{P}_k$, or $d(a, b) < d(x, y)$, then $A \cup B \in \mathcal{P}_{k+1}$. In other words, when all distances between members of $A$ and $B$ are smaller than all distances between points of other clusters, then $A$ and $B$ are joined next. Now, a cluster tree for $X$ is said to have a *cluster representation* $f : X \to \underline{Y}$ when all clustering methods that fulfil $(\star)$ reproduce this cluster tree when they are applied to the metric space $\underline{X}' := (X, d')$ with $d'(x, y) := e(fx, fy)$.

**Proposition 10** *Every cluster tree $\mathcal{P}_1, \ldots, \mathcal{P}_n$ for a finite set $X$ has a cluster representation in $\{0, \ldots, \lfloor (1+\sqrt{2})^n/4 \rfloor\}$ with Euclidean distance.*

*Proof.* Inductively, we construct maps $f_i : X \to \mathbb{Z}$ and integers $\delta_i$ such that $f_n$ is a cluster representation, and each $f_i$ is already "correct" for all $C \in \mathcal{P}_i$. For $C \in \mathcal{P}_i$, the convex hull of $f_i[C]$ will be the interval $[0, w_i(C)]$. For $A, B \in \mathcal{P}_i$ and $A \cup B \in \mathcal{P}_{i+1}$, $f_{i+1}[A \cup B]$ will be constructed by placing $f_i[A]$ and $f_i[B]$ besides each other at a distance $\delta_i$ that is larger than the diameter of any $C \in \mathcal{P}_i$, that is, with $\delta_i > w_i(C)$.

We start with $f_1(a) := 0$ for all $a \in X$, so that $w_1(A) = 0$ for all $A \in \mathcal{P}_1$, and put $\delta_1 := 1$. For $i \geqslant 1$, let $A_i, B_i \in \mathcal{P}_i$ be those elements with $C_i := A_i \cup B_i \in \mathcal{P}_{i+1}$ and $\min A_i < \min B_i$. Now put

$$\begin{aligned}
f_{i+1}(a) &:= f_i(a) \quad \text{for all } a \in A_i, \\
f_{i+1}(b) &:= f_i(b) + \delta_i + w_i(A) \quad \text{for all } b \in B_i, \\
f_{i+1}(x) &:= f_i(x) \quad \text{for all } x \notin C_i,
\end{aligned}$$

and $\delta_{i+1} := w_{i+1}(C_i) + 1$, where, by construction, $w_{i+1}(C_i) = \delta_i + w_i(A_i) + w_i(B_i)$. Then the convex hull of $f_{i+1}[C_i]$ is $[0, w_{i+1}(C_i)]$ as proposed. For all $C \in \mathcal{P}_{i+1}$ different from $C_i$, we have $C \in \mathcal{P}_i$ and thus $\delta_{i+1} > \delta_i > w_i(C) = w_{i+1}(C)$ as required. In case that $i \geqslant 2$, one of $A_i, B_i$ is in $\mathcal{P}_{i-1}$, hence either $w_i(A_i) = w_{i-1}(A_i)$ or $w_i(B_i) = w_{i-1}(B_i)$. Putting $m_i := \max\{w_i(A) : A \in \mathcal{P}_i\}$, this gives $m_{i+1} \leqslant 2m_i + m_{i-1} + 1$. It is easy to verify that the corresponding recursive upper bound $b_i$ with $b_{i+1} = 2b_i + b_{i-1} + 1$ and initial conditions $b_1 = 0$ and $b_2 = 1$ is $b_i = ((1+\sqrt{2})^i + (1-\sqrt{2})^i)/4 - 1/2 = \lfloor (1+\sqrt{2})^i/4 \rfloor$. In particular, $w_n(X) = m_n \leqslant b_n = \lfloor (1+\sqrt{2})^n/4 \rfloor$.

Finally, $f_n$ is a cluster representation: let $i \leqslant n$, $a \in A_i$, $b \in B_i$, $A' \neq B' \in \mathcal{P}_i$ with $\{A', B'\} \neq \{A_i, B_i\}$, and $a' \in A'$, $b' \in B'$. Then the smallest index $j$ for which there is $C \in \mathcal{P}_j$ with $a', b' \in C$ is at least $i+1$, hence $df_n(a,b) = df_i(a,b) < \delta_i \leqslant \delta_{j-1} \leqslant df_j(a',b') = df_n(a',b')$. □

Finally, this construction can be used to show that the following lower bound on order accuracy for maps into the real line:

**Theorem 11** *For every $n$-element metric space $\underline{X}$ with $n = 2^p$ for some integer $p$, there is a map $f : \underline{X} \to \mathbb{E}_1$ with order accuracy at least $3/7 - O(1/n)$.*

*Proof.* We iteratively define a binary cluster tree. For $k < n$, $\mathcal{P}_k$ is constructed from $\mathcal{P}_{k+1}$ as follows: choose some $C \in \mathcal{P}_{k+1}$ of maximal size, and let $w_C(\{x,y\})$ be the number of pairs $\{z,w\} \subseteq C$ with $0 < d(z,w) < d(x,y)$. In [PT86] it was proved that there is a partition of $C$ into two sets $A$ and $B$ of equal size such that

$$\sum_{x \in A,\, y \in B} w_C(\{x,y\}) \geqslant \frac{1}{2} \cdot \sum_{\{x,y\} \subseteq C} w_C(\{x,y\}) = \frac{1}{2} \cdot \binom{\binom{|C|}{2}}{2}.$$

Let $\mathcal{P}_k := \mathcal{P}_{k+1} \setminus \{C\} \cup \{A, B\}$. Note that $w_C(\{x,y\})$ is now the sum of $w_{A,B}(\{x,y\})$, the number of pairs $\{z,w\} \subseteq C$ with $0 < d(z,w) < d(x,y)$, $z \in A$, and $w \in B$, and of $w'_{A,B}(\{x,y\})$, the number of pairs $\{z,w\} \subseteq C$ with $0 < d(z,w) < d(x,y)$ and either $z, w \in A$ or $z, w \in B$.

Now we construct a representation as in the previous proposition, except that we might sometimes use $f'_i(a) := w_i(A_i) - f_i(a)$ and $f'_i(b) := w_i(B_i) - f_i(b)$ instead of $f_i(a)$ and $f_i(b)$ for the definition of $f_{i+1}|_{C_i}$. More precisely, when $f_i$ has already been defined and $A_i$, $B_i$, $C_i$ are as in the proposition, let $\gamma$ be the number of quadruples $(x, y, z, w) \in A_i \times B_i \times A_i \times B_i$ with $0 < d(z, w) < d(x, y)$ and $f_i(w) - f_i(z) < f_i(y) - f_i(x)$, and let $\gamma'$ be the number of quadruples $(x, y, z, w) \in A_i \times B_i \times A_i \times B_i$ with $0 < d(z, w) < d(x, y)$ and $f_i(z) - f_i(w) < f_i(x) - f_i(y)$. These numbers tell how many pairs of edges between $A_i$ and $B_i$ will be represented with the correct order of lengths when either $f_i$ or $f'_i$ is used for the definition of $f_{i+1}|_{C_i}$. Now put $f_{i+1}(x) := f_i(x)$ for all $x \notin C_i$, and either

$$\begin{aligned} f_{i+1}(a) &:= f_i(a) \quad \text{for all } a \in A_i, \text{ and} \\ f_{i+1}(b) &:= f_i(b) + \delta_i + w_i(A) \quad \text{for all } b \in B_i \end{aligned}$$

if $\gamma \geqslant \gamma'$, or otherwise

$$\begin{aligned} f_{i+1}(a) &:= f'_i(a) \quad \text{for all } a \in A_i, \text{ and} \\ f_{i+1}(b) &:= f'_i(b) + \delta_i + w_i(A) \quad \text{for all } b \in B_i. \end{aligned}$$

This assures that $|f_{i+1}(x) - f_{i+1}(y)| > |f_{i+1}(z) - f_{i+1}(w)|$ whenever $x \in A_i$, $y \in B_i$, and either $z, w \in A_i$ or $z, w \in B_i$. Moreover, since the sum of $\gamma$ and $\gamma'$ is $\binom{|A_i||B_i|}{2}$, their maximum is at least $|A_i||B_i|(|A_i||B_i| - 1)/4$. Hence, this step $i$ of the construction contributes to the overall accuracy $\alpha$ a summand $\alpha_i$ with

$$\alpha_i \cdot \binom{\binom{n}{2}}{2} \geqslant \sum_{x \in A_i,\, y \in B_i} w'_{A_i, B_i}(\{x, y\}) + \frac{|A_i||B_i|(|A_i||B_i| - 1)}{4}$$

$$= \sum_{x \in A_i,\, y \in B_i} \left( w_{C_i}(\{x, y\}) - w_{A_i, B_i}(\{x, y\}) \right) + \frac{|A_i||B_i|(|A_i||B_i| - 1)}{4}$$

$$= \sum_{x \in A_i,\, y \in B_i} w_{C_i}(\{x, y\}) - \binom{|A_i||B_i|}{2} + \frac{|A_i||B_i|(|A_i||B_i| - 1)}{4}$$

$$\geqslant \frac{1}{2} \cdot \binom{\binom{|C_i|}{2}}{2} - \frac{|A_i||B_i|(|A_i||B_i| - 1)}{4} = \frac{3}{64}|C_i|^4 + O(|C_i|^3).$$

Finally, all $C_i$ are of size $n/2^q$ for some $q$ with $0 \leqslant q < p$, and there are exactly $2^q$ many of this size. Hence the overall accuracy is

$$\alpha = \sum_{i=1}^{n} \alpha_i \geqslant \sum_{q=0}^{p-1} 2^q \cdot \frac{3}{8}(1/2^q)^4 + O(1/n) = \frac{3}{7} - O(1/n).$$

$\square$

However, this lower bound is very likely not the best possible. The rank correlation $\varrho$ between two independently chosen linear orders on $m$ elements is nearly normally distributed with expected value 0 and standard deviation $O(1/\sqrt{m})$ (cf. [KG90]). Hence $(\varrho + 1)/2$ has expected value $1/2$, which motivates the following conjecture.

**Conjecture 12** *Every finite metric space can be mapped into $\mathbb{E}_1$ with accuracy $\geqslant 1/2$.*